\documentclass[11pt]{amsart}

\usepackage{amsmath}
\usepackage{amssymb}
\usepackage{mathrsfs}
\usepackage{comment}
\usepackage{color}
\usepackage[colorlinks,citecolor=blue,urlcolor=black,linkcolor=black]{hyperref}

\newtheorem{theorem}{Theorem}[section]

\newtheorem{lemma}[theorem]{Lemma}

\newtheorem{proposition}[theorem]{Proposition}
\newtheorem{corollary}[theorem]{Corollary}

\theoremstyle{definition}
\newtheorem{definition}[theorem]{Definition}

\newtheorem{question}[theorem]{Question}

\theoremstyle{remark}
\newtheorem{remark}[theorem]{Remark}

\newcount\skewfactor
\def\mathunderaccent#1#2 {\let\theaccent#1\skewfactor#2
\mathpalette\putaccentunder}
\def\putaccentunder#1#2{\oalign{$#1#2$\crcr\hidewidth
\vbox to.2ex{\hbox{$#1\skew\skewfactor\theaccent{}$}\vss}\hidewidth}}
\def\name{\mathunderaccent\tilde-3 }

\def\smallbox#1{\leavevmode\thinspace\hbox{\vrule\vtop{\vbox
   {\hrule\kern1pt\hbox{\vphantom{\tt/}\thinspace{\tt#1}\thinspace}}
   \kern1pt\hrule}\vrule}\thinspace}

\newcommand{\cf}{{\rm cf}}
\newcommand{\ra}{\rangle}
\newcommand{\la}{\langle}

\def\qedref#1{$\qed_{\reforiginal{#1}}$}

\title{J\'onsson and Magidor filters}
\author{Omer Ben-Neria}
\address{Einstein Institute of Mathematics
 The Hebrew University of Jerusalem,
 Jerusalem 91904, Israel}
\email{omer.bn@mail.huji.ac.il}
\thanks{Ben-Neria would like to thank the Israel Science Foundation (Grants 1832/19 and 1302/23) for
their support}

\author{Shimon Garti}
\address{Einstein Institute of Mathematics,
 The Hebrew University of Jerusalem,
 Jerusalem 91904, Israel}
\email{shimon.garty@mail.huji.ac.il}
\thanks{Garti would like to thank the Israel Science Foundation (Grant 2320/23) for
their support}


\begin{document}
\let\labeloriginal\label
\let\reforiginal\ref
\def\ref#1{\reforiginal{#1}}
\def\label#1{\labeloriginal{#1}}

\begin{abstract}
We study the J\'{o}nsson and Magidor partition properties and their strengthened filter forms, primarily in the context of the Axiom of Determinacy.
By applying methods of Kleinberg concerning strong partition cardinals $\kappa$ and their associated cardinal sequence $(\kappa_n : n \in \omega)$, we show that the limit cardinal $\kappa_\omega$ is a Magidor cardinal if $\kappa > \aleph_1$.
We also force over a model of \textsf{AD} and $V=L(\mathbb{R})$ to get the existence of a singular cardinal $\lambda$ carrying a Magidor filter, addressing a question from \cite{MR3666820}.
\end{abstract}

\maketitle

\newpage

\section{Introduction}

In this work, we consider strong partition properties of J\'{o}nsson and Magidor, and the existence of associated filters, both under the axiom of choice and in the choiceless context, under the axiom of determinacy.

An infinite cardinal $\lambda$ is J\'onsson iff $\lambda\rightarrow[\lambda]^{<\omega}_\lambda$, that is for every coloring $c:[\lambda]^{<\omega}\rightarrow\lambda$ one can find $A\in[\lambda]^\lambda$ such that $c''[A]^{<\omega}\neq\lambda$, in which case we shall say that $A$ is $c$-monochromatic.
An infinite cardinal $\lambda$ is $\omega$-J\'onsson iff $\lambda\rightarrow[\lambda]^{\omega}_\lambda$.
Erd\H{o}s and Hajnal proved in \cite{MR0209161} that there are no $\omega$-J\'onsson cardinals in \textsf{ZFC}.
In the absence of the axiom of choice, for example under \textsf{AD}, the stronger notion of $\omega$-J\'onsson cardinals may be consistent.

An intermediate concept between the notions of J\'{o}nsson and $\omega$-J\'{o}nsson cardinals is the one of Magidor cardinal (Magidority), which was developed in \cite{MR3750266}.
A cardinal $\lambda$ is Magidor iff $\lambda\rightarrow[\lambda]^{\aleph_0\text{-bd}}_\lambda$, where $\aleph_0\text{-bd}$ refers to all bounded subsets of $\lambda$ of cardinality $\aleph_0$.
Many results concerning the possible strength and obstructions to Magidor cardinals appear when considering the ordertype restricted notion of $\omega$-Magidority.  We say that $\lambda$ is $\omega$-Magidor iff $\lambda\rightarrow[\lambda]^{\omega\text{-bd}}_\lambda$
where $\omega\text{-bd}$ refers to all bounded subsets of $\lambda$ of ordertype $\omega$.
Magidor observed that if $\lambda$ satisfies the large cardinal property I1 then $\lambda$ is Magidor, see \cite[Question 24.1]{MR1994835}, hence Magidority is consistent with \textsf{ZFC}.

It is known that under \textsf{AD}, many cardinals below $\Theta$ satisfy many partition properties.
Considering J\'{o}nssonicity,
Kleinberg proved that the first uncountable cardinals up to and including $\aleph_\omega$ are J\'{o}nsson (\cite{MR1994835}). Jackson extended the result to all cardinals $\kappa < \aleph_{\omega_1}$,
Steel (see \cite{MR2768698}) proved that every regular uncountable cardinal below $\Theta$ is measurable and thus J\'{o}nsson.
Jackson, Ketchersid, Schlutzenberg, and Woodin, have independently proved that the results extend to all uncountable cardinals below $\Theta$ (\cite{MR3343535}).
In this work we incorporate Magidority into the analysis of large cardinals under determinacy.
Strong results in \textsf{ZF} concerning $\omega$-J\'{o}nsson cardinals appear in \cite{MR2135668}.

We consider here the question of whether \textsf{AD} implies that certain singular cardinals below $\Theta$ are Magidor.
We show that under \textsf{AD}, $\aleph_\omega$ is $\omega$-Magidor and that for every cardinal $\kappa$, satisfying $\aleph_1 < \kappa < \Theta$ and
$\kappa\rightarrow (\kappa)^\kappa_2$, its associated singular cardinal $\kappa_\omega$ is Magidor. In particular $\aleph_{\omega+\omega}$ which is associated to $\kappa = \aleph_{\omega+1}$ is Magidor.
Our results in the AD context build on  the work of Kleinberg (\cite{MR1994835}) on cardinals $\kappa < \Theta$ satisfying the strong partition property, and their associated ultrapowers sequence $(\kappa_n : n < \omega)$.

A natural strengthening of square bracket relations is given by the following filter version of these notions.

\begin{definition}
\label{deffilters}{(Square brackets filters)} \newline
Let $\mathscr{F}$ be a uniform filter over $\lambda$ which contains all the end-segments of $\lambda$.
\begin{enumerate}
\item [$(\aleph)$] We shall say that $\mathscr{F}$ is J\'onsson iff for every $c:[\lambda]^{<\omega}\rightarrow\lambda$ one can find $A\in\mathscr{F}$ such that $c''[A]^{<\omega}\neq\lambda$.
\item [$(\beth)$] We shall say that $\mathscr{F}$ is Magidor iff for every $c:[\lambda]^{\aleph_0\text{-bd}}\rightarrow\lambda$ one can find $A\in\mathscr{F}$ such that $c''[A]^{\aleph_0\text{-bd}}\neq\lambda$.
\end{enumerate}
\end{definition}

The filter version of square brackets relations is usually stronger than the classical one.
In the particular case of Magidority, it is known that Magidor filters do not exist in \textsf{ZFC} despite the fact that there might be Magidor cardinals.
However, we shall force the existence of Magidor filters over a model of \textsf{AD}.

The paper is organized as follows:\\

\noindent
In \textbf{Section \ref{sec:JonssonFiltersAC}}, we consider the notion of J\'{o}nsson filters $\mathscr{F}$ and their J\'{o}nsson degrees $\alpha_J(\mathscr{F})$ (see Lemma \ref{aalphafilter}) in the ZFC context. We show that a singular limit of cardinals carrying J\'{o}nsson filters of bounded degree must also carry a J\'{o}nsson filter. This generalizes a well-known theorem of Prikry that a singular limit of measurable cardinals carries a J\'{o}nsson filter.\\

\noindent
In \textbf{Section \ref{sec:MagidorSingularAD}} we extend the analysis of Kleinberg sequences $\la \kappa_n : n < \omega\ra$ with respect to Magidority, and show that if $\kappa > \aleph_1$ is a strong partition cardinal, then the associated limit of its Kleinberg sequence $\kappa_\omega$ is Magidor, and that $\aleph_\omega$ is $\omega$-Magidor.
We also show that each $\kappa_n$, $2 < n < \omega$ in the Kleinberg sequence is not Magidor.
\\

\noindent
Finally, in \textbf{Section \ref{sec:ForcingMagidorFilters}} we consider the Magidor filter property at Prikry generic extensions of \textsf{AD} models. This follows a long sequence of combinatorial and arithmetical results obtained by forcing with Prikry-like posets over models of \textsf{AD}:  Apter \cite{MR604883}, Spector \cite{MR583378}, Henle \cite{MR722169}, and Apter-Henle-Jackson \cite{MR1695015}.
We show that a singular cardinal carrying a Magidor filter may exist in such extensions, providing a partial answer to \cite[Question 2.5]{MR3666820} concerning the existence of Magidor filters in models of \textsf{AD}.

\newpage

\section{Sufficient conditions for J\'{o}nsson filters over singular cardinals}\label{sec:JonssonFiltersAC}

In this section we work in \textsf{ZFC} and we follow in the footsteps of Prikry (\cite{MR0262075}) in order to show that a singular cardinal which is a limit of cardinals carrying J\'{o}nsson filters (with some parameter) must also carry a J\'{o}nsson filter. In the next section we shall prove a similar theorem under \textsf{AD}.

Prikry used normal ultrafilters over the measurable cardinals along the way, but the proof needs only the partition theorem of Rowbottom for such ultrafilters. In particular, one can replace the normal ultrafilters by Rowbottom or J\'{o}nsson filters. This shift is quite meaningful since the existence of a normal ultrafilter is expensive in terms of consistency strength.
The proof below follows the ideas of Prikry \cite{MR0262075} and the presentation of Kanamori \cite{MR1994835}.

We phrase the theorem (and the proof) in the more general context of $\theta$-Rowbottom cardinals.
Rowbottom cardinals possess two equivalent definitions. The first one is model theoretic. We say that $\kappa$ is Rowbottom iff for every structure $\mathfrak{A}=\langle A,Q,\ldots\rangle$ for a countable language such that $|A|=\kappa$, $|Q| < \kappa$, there is an elementary substructure $\mathfrak{B}=\langle B,Q\cap B,\ldots\rangle$ so that $|B|=\kappa$ and $|Q\cap B|\leq\aleph_0$.
One may add another parameter, and define $\kappa$ as a $\nu$-Rowbottom cardinal exactly in the same way but with the demand $|B|=\kappa$ and $|Q\cap B|<\nu$.

The second natural definition of Rowbottom cardinals is combinatorial, and based on the square brackets notation. We say that $\kappa$ is Rowbottom iff $\kappa\rightarrow [\kappa]^{<\omega}_{\gamma,<\omega_1}$ for all $\gamma<\kappa$. This means that for every $\gamma<\kappa$ and a coloring $f:[\kappa]^{<\omega}\rightarrow \gamma$, there exists a set $A\in[\kappa]^\kappa$ for which $f''[A]^{<\omega}$ is countable.
Again, we can generalize the definition into $\nu$-Rowbottomness upon replacing $\omega_1$ by $\nu$, i.e. $\kappa\rightarrow [\kappa]^{<\omega}_{\gamma,<\nu}$.
We note that if $\kappa>\cf(\kappa)$ is such a cardinal then $\cf(\kappa)<\nu$ (\cite{MR1994835}).
In particular, if  $\kappa$ is a singular Rowbottom cardinal then $\kappa$ has countable cofinality.

\begin{theorem}
\label{mtrowbottom} (Rowbottom filters over singular cardinals) \newline
Assume that:
\begin{enumerate}
\item [$(\aleph)$] $\cf(\kappa) = \nu\leq\theta<\kappa$.
\item [$(\beth)$] $\langle\kappa_\alpha:\alpha<\nu\rangle$ is a continuous increasing sequence of cardinals which tends to $\kappa$.
\item [$(\gimel)$] $\mathscr{F}_{\alpha+1}$ is a $\kappa_{\alpha+1}$-complete $\theta^+$-Rowbottom filter over $\kappa_{\alpha+1}$ for every $\alpha<\nu$.
\item [$(\daleth)$] $\gamma<\kappa\Rightarrow \gamma^\theta<\kappa$.
\end{enumerate}
Then there exists a $\theta^+$-Rowbottom filter over $\kappa$.
\end{theorem}

\par\noindent\emph{Proof}. \newline
We define a uniform filter $D$ over $\kappa$ as follows:

$$
x\in D\Leftrightarrow \exists\beta\forall\alpha\in[\beta,\nu), x\cap\kappa_{\alpha+1}\in \mathscr{F}_{\alpha+1}.
$$

Our objective is to show that $D$ is a $\theta^+$-Rowbottom filter.
The degree of completeness of $D$ is low, and in the typical case of Rowbottom cardinals (i.e. $\cf(\kappa)=\omega$) it is not even $\aleph_1$-complete. This fact is problematic, since the completeness of the filter is used in most proofs of constructing monochromatic sets.
The solution of this problem is to choose our large sets in such a way that the ordinal $\beta$ in the above definition of $D$ is constant.\\

We begin by the following definition of type. For every $s\in[\kappa]^{<\omega}$ let ${\rm type}(s)$ be the set $\{(\alpha,\ell)\in\nu\times\omega: 0<\ell=|s\cap (\kappa_{\alpha+1}-\kappa_\alpha)|\}$. Assume $\gamma<\kappa$ and $f:[\kappa]^{<\omega}\rightarrow\gamma$ is any function. Let $T$ be a type (i.e., a set of ordered pairs from $\nu\times\omega$), let $H$ be a subset of $\kappa$ and let $W^f_H(T)$ be the set $\{\zeta<\gamma:\exists s\in[H]^{<\omega}, T={\rm type}(s), f(s)=\zeta\}$. We shall try to prove the following:

\begin{center}
For every $\gamma<\kappa$ and every $f:[\kappa]^{<\omega}\rightarrow\gamma$ there exists $H\in D$ \\ such that $|W^f_H(T)|\leq\theta$ for every type $T$.
\end{center}

If we succeed then the proof of the theorem will be accomplished. Indeed, the number of types is $|\nu\times\omega|$, and for every type there are at most $\theta$ values for $f$ on $H$, so $|f''[H]^{<\omega}| \leq\nu\cdot\theta=\theta\cdot\theta<\kappa$. Since $H\in D$ we conclude that $D$ is a $\theta^+$-Rowbottom filter.

For proving the above claim, we decompose it into $\aleph_0$ local assertions, which we shall prove by induction. Let $(\star)_n$ be the assertion that for every $\gamma<\kappa$ and every $f:[\kappa]^{<\omega}\rightarrow\gamma$ there is a subset $H_n\subseteq\kappa$ such that:
\begin{enumerate}
\item [$(\aleph)_n$] $H_n\cap\kappa_{\alpha+1}\in \mathscr{F}_{\alpha+1}$ whenever $\kappa_{\alpha+1}>\gamma$.
\item [$(\beth)_n$] If $T$ is any type such that $|T|=n$ then $|W^f_{H_n}(T)|\leq\theta$.
\end{enumerate}
If we prove $(\star)_n$ for every $n\in\omega$ then we are done. We will choose, for every $n\in\omega$, a set $H_n\subseteq\kappa$ which satisfies $(\star)_n$. Then we will set $H=\bigcap\{H_n:n\in\omega\}$. Although $D$ is not $\aleph_1$-complete in general, here we have the same end-segment for every $H_n$ (it is determined by the fixed number of colors, $\gamma$). The completeness of each $\mathscr{F}_{\alpha+1}$ and item $(\aleph)_n$ for every $n\in\omega$ make sure that $H\in D$. From $(\beth)_n$ for every $n\in\omega$ it follows that $|W^f_H(T)|\leq\theta$ for every type $T$, since there must be some $n$ which is the type of $T$.

The proof of $(\star)_n$ is by induction on $n\in\omega$. If $n=0$ then $T=\emptyset$ and hence only $s=\emptyset$ is of type $T$. We can choose $H_0=\kappa\in D$ and the requirements hold trivially. So we focus on $(\star)_{n+1}$, assuming that $(\star)_n$ has been shown.

Recall that $f:[\kappa]^{<\omega}\rightarrow\gamma$ for some $\gamma<\kappa$. For every $\alpha<\nu$ such that $\gamma<\kappa_{\alpha+1}$ and every $t\in[\kappa_\alpha]^{<\omega}$ we define a function $f_t:[\kappa_{\alpha+1}- \kappa_\alpha]^{<\omega} \rightarrow\gamma$ by $f_t(u)=f(t\cup u)$. Assumption $(\gimel)$ in the statement of the theorem provides us with $x^\alpha_t\in \mathscr{F}_{\alpha+1}$ such that $|f''_t[x^\alpha_t]^{<\omega}|\leq\theta$. We let $x_\alpha=\bigcap\{x^\alpha_t: t\in[\kappa_\alpha]^{<\omega}\}-\kappa_\alpha$. Notice that $x_\alpha\in \mathscr{F}_{\alpha+1}$ for every pertinent $\alpha<\nu$.

By $(\daleth)$, we may assume that $\gamma^\theta=\gamma$.
Fix a bijection $h:[\gamma]^{\leq\theta}\rightarrow\gamma$. For every $(\alpha,\ell)\in\nu\times\omega$ we define $f_{\alpha\ell}:[\kappa_\alpha]^{<\omega} \rightarrow[\gamma]^{\leq\theta}$ by
\[
f_{\alpha\ell}(t)= f_t''[x_\alpha]^{<\omega} := \{\eta<\gamma \mid \exists u \in [x_\alpha]^{<\omega}, f_t(u)=\eta\}.
\]
We convert $f_{\alpha\ell}$ which is a set-mapping into a function $g_{\alpha\ell}:[\kappa_\alpha]^{<\omega} \rightarrow\gamma$ by letting $g_{\alpha\ell}(t)=h(f_{\alpha\ell}(t))$.

For every $(\alpha,\ell)\in\nu\times\omega$ we use the induction hypothesis $(\star)_n$ with respect to $g_{\alpha\ell}$ in order to elicit $H_{\alpha\ell}\subseteq\kappa$ such that:
\begin{enumerate}
\item [$(a)$] $H_{\alpha\ell}\cap\kappa_{\alpha+1}\in \mathscr{F}_{\alpha+1}$ whenever $\gamma<\kappa_{\alpha+1}$.
\item [$(b)$] If $T'$ is a type of size $n$ then $|W^{g_{\alpha\ell}}_{H_{\alpha\ell}}(T')|\leq\theta$.
\end{enumerate}
Finally, we can define the following set:

$$
H_{n+1}=\bigcup\{x_\alpha:\gamma<\kappa_{\alpha+1}\}\cap\bigcap\{H_{\alpha\ell}: (\alpha,\ell)\in\nu\times\omega\}.
$$

First we observe that $H_{n+1}\in D$. Indeed, $H_{n+1}\cap\kappa_{\alpha+1}\in \mathscr{F}_{\alpha+1}$ whenever $\gamma<\kappa_{\alpha+1}$ since $x_\alpha\in \mathscr{F}_{\alpha+1}$ and every $H_{\beta\ell}$ belongs to $\mathscr{F}_{\alpha+1}$. Consequently, $H_{n+1}\in D$ as required in the first part of $(\star)_{n+1}$.

Second, assume that $T$ is any type of size $n+1$ and let $\alpha$ be the maximal ordinal for which $(\alpha,\ell)\in T$. If ${\rm type}(s)=T$ then $s=t\cup u$ for $t,u$ such that $|{\rm type}(t)|=n$ and ${\rm type}(u)=\{(\alpha,\ell)\}$.
Let $T'$ be ${\rm type}(t)$.
By the induction hypothesis applied to $T'$ we have $|W^{g_{\alpha\ell}}_{H_{\alpha\ell}}(T')|\leq\theta$.

Now $f(s)=f(t\cup u)=f_t(u)\in f_{\alpha\ell}(t)$. Recall that $h(f_{\alpha\ell}(t))=g_{\alpha\ell}(t)\in W^{g_{\alpha\ell}}_{H_{\alpha\ell}}(T')$. Since $|W^{g_{\alpha\ell}}_{H_{\alpha\ell}}(T')|\leq\theta$ and $|f_{\alpha\ell}(t)|\leq\theta$ we infer that $f(s)$ assumes at most $\theta\cdot\theta$ possibilities and hence $|W_{H_{n+1}}^f(T)|\leq\theta$ as desired.

\hfill \qedref{mtrowbottom}

By \cite{MR654576}, if $\kappa=\cf(\kappa)>2^\omega$ is a Rowbottom cardinal then $\kappa$ is strongly inaccessible, and hence item $(\daleth)$ in the above theorem is satisfied for these cardinals. We may conclude:

\begin{corollary}
\label{cor1} Assume $\kappa>2^\omega$ is a singular cardinal of countable cofinality, a limit of regular cardinals each of which carries a Rowbottom filter. Then there is a Rowbottom filter over $\kappa$ as well.
\end{corollary}

\hfill \qedref{cor1}

\begin{remark}
\label{rrr} If $\mathscr{U}$ is a normal ultrafilter over $\kappa$ then, in $L[\mathscr{U}]$, every J\'onsson cardinal is regular, as proved in \cite{MR0277346}.
It follows from this fact and the above theorem that the number of cardinals which carry a Rowbottom filter in $L[\mathscr{U}]$ is finite. However, this is subsumed by the known fact that $\kappa$ (being measurable) is the only cardinal which carries a Rowbottom or J\'onsson filter in $L[\mathscr{U}]$.
\end{remark}

\hfill \qedref{rrr}

We turn to J\'{o}nsson filters.
We would like to prove a similar theorem about J\'{o}nsson filters at singular cardinals. We need a parallel to the parameter $\theta$ of the Rowbottomness degree of the filters. This parameter is given by the following version of an argument by Kleinberg (see \cite{MR1994835}):

\begin{lemma}
\label{aalphafilter} Assume $F$ is a J\'onsson filter over $\kappa$. \newline
Then there exists a cardinal $\alpha<\kappa$ such that for every $f:[\kappa]^{<\omega}\rightarrow\alpha$ there is an element $x\in F$ for which $|f''[x]^{<\omega}|<\alpha$.
\end{lemma}

\par\noindent\emph{Proof}. \newline
We begin by showing that there exists an ordinal $\alpha<\kappa$ such that for every $f:[\kappa]^{<\omega}\rightarrow\alpha$ there is an element $x\in F$ for which $f''[x]^{<\omega}\neq\alpha$. Towards contradiction assume that for every $\alpha<\kappa$ there is a function $f_\alpha:[\kappa]^{<\omega} \rightarrow\alpha$ such that $f''_\alpha[x]^{<\omega}=\alpha$ for every $x\in F$. Define a function $f:[\kappa]^{<\omega}\rightarrow\kappa$ by $f(\emptyset)=0$ and if $t\neq\emptyset, \alpha={\rm min}(t)$ then $f(t)=f_\alpha(t-\{\alpha\})$.

Choose any element $x\in F$ and any color $\gamma<\kappa$. Pick an ordinal $\gamma<\alpha\in x$ and let $y=x-(\alpha+1)$. By the assumptions on $F$ we know that $y\in F$ and hence there is $s\in[y]^{<\omega}$ such that $f_\alpha(s)=\gamma$. Let $t=s\cup\{\alpha\}$, so $\alpha={\rm min}(t)$. It follows that $f(t)=f_\alpha(t-\{\alpha\})=f_\alpha(s)=\gamma$. However, $t=s\cup\{\alpha\}\subseteq x$ and hence $f''[x]^{<\omega}=\kappa$ for every $x\in F$, a contradiction.

Now let $\alpha<\kappa$ be the first ordinal with the above property. It follows that $\alpha$ must be a cardinal (if not, use a bijection from some $\beta<\alpha$ in order to show that $\alpha$ is not the first with the above property) and we shall prove that $\alpha$ satisfies the lemma. By the minimaity assumption of $\alpha$, we can choose for every $\beta<\alpha$ a function $f_\beta:[\kappa]^{<\omega}\rightarrow\beta$ so that $f''_\beta[x]^{<\omega}=\beta$ whenever $x\in F$.

Towards contradiction assume that there exists $g:[\kappa]^{<\omega} \rightarrow\alpha$ such that $|g''[x]^{<\omega}|=\alpha$ for every $x\in F$. We wish to define a new function $f:[\kappa]^{<\omega} \rightarrow\alpha$, which omits no color on every $y\in F$.
Given $t=\{t_1,\ldots,t_n\}\in[\kappa]^{<\omega}$, if $n$ is not of the form $2^i\cdot 3^j$ then $f(t)=0$. If $n=2^i\cdot 3^j$ then $f(t)=f_{g(t_1,\ldots,t_i)} (t_{i+1},\ldots,t_{i+j})$.

Let $y$ be any element of $F$, and $\gamma<\alpha$ a color. Choose $\{t_1,\ldots,t_i\}\in [y]^{<\omega}$ and $\beta\in(\gamma,\alpha)$ for which $g(\{t_1,\ldots,t_i\})=\beta$. By the above notation, $\gamma\in f''_\beta[y]^{<\omega}$. Hence, there is $s=\{s_1,\ldots,s_j\}\subseteq y$ such that $f_\beta(s)=\gamma$. Without loss of generality, $t_i<s_1$. Define $r=t\cup s$ and observe that $f(r)=f_\beta(s)=\gamma$, a contradiction.

\hfill \qedref{aalphafilter}

For a J\'onsson filter $\mathscr{F}$ let $\alpha_J(\mathscr{F})$ be the first ordinal $\alpha$ whose existence is proved by the above lemma.
The parallel to Theorem \ref{mtrowbottom} follows, by assuming that all the filters share the same $\alpha_J(\mathscr{F})$. We phrase it as a corollary:

\begin{corollary}
\label{corforjfilters} J\'{o}nsson filters over singular cardinals. \newline
Assume that:
\begin{enumerate}
\item [$(\aleph)$] $\nu=\cf(\kappa)<\kappa$.
\item [$(\beth)$] $\langle\kappa_\alpha:\alpha<\nu\rangle$ is a continuous increasing sequence of cardinals which tends to $\kappa$.
\item [$(\gimel)$] $\mathscr{F}_{\alpha+1}$ is a $\kappa_{\alpha+1}$-complete J\'onsson filter over $\kappa_{\alpha+1}$ which contains the end-segments, for every $\alpha<\nu$.
\item [$(\daleth)$] $\alpha_J(\mathscr{F}_{\alpha+1})=\theta$ for some fixed $\theta\in [\nu,\kappa)$, for every $\alpha<\nu$.
\end{enumerate}
Then there exists a J\'{o}nsson filter over $\kappa$.
\end{corollary}

\hfill \qedref{corforjfilters}

The above corollary demonstrates the similarity between Rowbottomness and J\'{o}nssonicity.
According to the formal definition, the colorings in J\'{o}nsson cardinals are from $[\lambda]^{<\omega}$ into $\lambda$.
such colorings may omit colors but one cannot expect a set which assumes \emph{less than} $\lambda$ colors.
To see this, consider the coloring defined by $c(\bar{\alpha})=\min(\bar{\alpha})$.
Rowbottomness, therefore, applies to colorings from $[\lambda]^{<\omega}$ into $\gamma$ where $\gamma<\lambda$.

However, if the range of the coloring is less than $\lambda$ then the phenomenon of reducing the number of colors to a set of size strictly less than the range of the coloring is manifested at J\'{o}nsson cardinals (or filters) as well.
We know that for some $\alpha\in\lambda$ one has $\lambda\rightarrow[\lambda]^{<\omega}_{\alpha,<\alpha}$.
The only point is that if $\lambda$ is Rowbottom then this will happen with a fixed size of colors obtained at every coloring into every $\gamma<\lambda$, while J\'{o}nssonicity is weaker from this point of view.

\begin{remark}
Corollary \ref{corforjfilters} does not extend to a regular limit of cardinals carrying  J\'{o}nsson filters or even to a regular limit of measurable cardinals.
Indeed, from a suitable large cardinal assumption, it is possible to force a regular cardinal $\kappa$ being a limit of measurable cardinals $\kappa=\bigcup_{\alpha\in\kappa}\kappa_\alpha$ and carrying a nonreflecting stationary set $S \subseteq \kappa$.
Thus, by a theorem of Tryba and Woodin, $\kappa$ is not J\'onsson.
\end{remark}

\newpage

\section{Magidority and singular cardinals}\label{sec:MagidorSingularAD}

In this section we extend the Kleinberg and Henle analysis to study the Magidor partition property.
Firstly we show that if $\kappa$ is a strong partition cardinal and $(\kappa_n:n\in\omega)$ is the corresponding Kleinberg sequence then $\kappa_n$
is not Magidor for every $n>2$.
We then show that the limit $\kappa_\omega$ of the Kleinberg sequence is a Magidor cardinal if $\kappa > \aleph_1$, and is $\omega$-Magidor if $\kappa = \aleph_1$ (i.e., $\kappa_\omega = \aleph_\omega$).
We commence with relevant background on the Kleinberg sequence, the information is based on \cite{MR0479903}.

For a subset $D \subseteq \kappa$ of size $\kappa$,
let $[D]^\kappa$ denote the set of all increasing functions $f : \kappa \to D$.
For a function $h \in [\kappa]^\kappa$ let  ${}_\omega h$ be the function which enumerates the $\omega$-club of all limits ordinal $\alpha \in \kappa$ which are limits of countable strictly increasing sequences $\alpha = \cup_{n\in\omega} h(\alpha_n)$, of ordinals $\alpha_n \in \kappa$.
Clearly, ${}_\omega h \in \mathscr{U}$.

A cardinal $\kappa > \aleph_0$ is a \emph{strong partition cardinal} iff $\kappa\rightarrow(\kappa)^\kappa$. Such cardinals are measurable, as witnessed by the restriction of the club filter to the set of countable cofinality ordinals $\alpha < \kappa$, which we denote by $\mathscr{U}$.
The filter $\mathscr{U}$ is a $\sigma$-complete normal measure on $\kappa$. Moreover, it is a strong partition (ultra)filter over $\kappa$ in the sense that for every coloring function $c : [\kappa]^\kappa \to \nu$ for some $\nu < \kappa$, there exists a some $h \in [\kappa]^\kappa$ so that the restriction of $c$ to $[{}_\omega h]$ is monochromatic.

One assigns to a strong partition cardinal $\kappa$ an increasing sequence of cardinals $(\kappa_n : n \in \omega)$. The informal idea of the definition is to start with $\kappa_1 = \kappa$ and for each $n \in \omega$, use ultrapowers of $\kappa_n$ by $\mathscr{U}$ to define  $\kappa_{n+1}$ to be the ordertype of the well-ordering of $\mathscr{U}$-equivalence classes of ${}^\kappa \kappa_n$.

The second cardinal in the sequence $\kappa_2$ can be shown to be regular and even measurable.
Moreover, it satisfies the weak partition property denoted by $\kappa_2\rightarrow(\kappa_2)^{<\kappa_2}$ (making $\kappa_2$ a \emph{weak partition cardinal}).
The other cardinals $\kappa_n$, $n \geq 3$ are shown to be singular of cofinality $\cf(\kappa_n) = \kappa_2$.\footnote{Moreover, by the work of Kleinberg and Henle these cardinals satisfy relatively strong square bracket partition relations. In particular, under \textsf{AD} each $\aleph_n$ is J\'onsson and in fact each $\aleph_n$ carries a J\'onsson filter as shown in \cite{MR0479903}. In this specific case, the Kleinberg sequence begins at $\aleph_1$.}

Kleinberg showed that many combinatorial properties of the $\kappa_n$s can be translated to properties of the strong partition cardinal $\kappa$.
Doing so requires working with representing functions for ordinals below $\kappa_n$, and may face difficulties involving choosing representing functions from the associated $\mathscr{U}$-equivalences classes.
To circumvent choice related difficulties, Kleinberg defines $\kappa_n$ directly from functions $F \in {}^{\kappa^{n-2}}\kappa$ and the natural $\mathscr{U}$-induced equivalence relation.

In this section we try to apply a similar method with respect to Magidority.
It turns out that the situation is not the same.
An interesting important difference is that the $\kappa_n$s in the Kleinberg sequence are not Magidor.
Ahead of proving this fact we shall need the following:

\begin{definition}\label{def:beta}${}$
\begin{enumerate}
    \item
    For $n \geq 3$ and $\gamma \in \kappa_n$, represented as $[F]_{\mathscr{U}}$ for some $F \in [\kappa]^{\kappa^{n-2}}$, we define for every
    $1 \leq k \leq n-3$ the
    induced $(-k)$-representation of $\gamma$, to be the function $\beta^F_{-k} \in [\kappa]^k \to \kappa_{n-k}$, defined byW
    \[
    \beta^{F}_{-k}(\alpha_0,\dots, \alpha_{k-1}) = [F_{\alpha_0,\dots,\alpha_{k-1}}]_{\mathscr{U}}
    \]
    where $F_{\alpha_0,\dots,\alpha_{k-1}} \in [\kappa]^{\kappa_{n-k}}$ is given by
    \[
    F_{\alpha_0,\dots,\alpha_{k-1}}(\alpha_k,\dots,\alpha_{n-1}) = F(\alpha_0,\dots,\alpha_{k-1},\alpha_k,\dots,\alpha_{n-1})
    \]

    \item Let $\gamma_0,\gamma_1 < \kappa_n$ and $1 \leq k < n-2$, where $F_0,F_1$ represent $\gamma_0,\gamma_1$ respectively. We say that $\gamma_0,\gamma_1$ are $k$-interlaced if
    \[
    \min_{A \in \mathscr{U}^k}(\sup_{(\alpha_0,\dots,\alpha_{k-1}) \in A}\beta^{F_0}_{-k}(\alpha_0,\dots,\alpha_{k-1})) =
    \min_{A \in \mathscr{U}^k}(\sup_{(\alpha_0,\dots,\alpha_{k-1}) \in A}\beta^{F_1}_{-k}(\alpha_0,\dots,\alpha_{k-1}))
    \]
Where $\mathscr{U}^k$ is the usual product measure of $\mathscr{U}$.
\noindent
The \textbf{interlacing number} of $\gamma_0,\gamma_1$ is the maximal $0 \leq k < n-2$ for which $\gamma_0,\gamma_1$ are $k$-interlaced.
It is straightforward to see that the definition does not depend on the choice of the representing functions.

\item Let $\gamma_0,\gamma_1 < \kappa_\omega$ be a pair of ordinals, with $n_0 \leq n_1 < \omega$ minimal so that $\gamma_i \in \kappa_{n_i}$, $i < 2$, and
$F_0,F_1$ representing functions $F_i \in [\kappa]^{\kappa_{n_i}}$, $i<2$, respectively.
If $n_0 < n_1$, let
\[A_{F_0,F_1} = \{ (\alpha_0,\dots, \alpha_{n_1-n_0-1}) \in [\kappa]^{n_1-n_0}:
 \gamma_0 < \beta^{F_1}_{n_1-n_0}(\alpha_0,\dots,\alpha_{n_1-n_0-1})\}\]

 and for each $k < n_0-2$,
 $B_{F_0,F_1,k} \subseteq [\kappa]^{n_1-n_0}$
 be the set of all $(\alpha_0,\dots, \alpha_{n_1-n_0-1}) \in [\kappa]^{n_1-n_0}$
 for which the interlacing number of
 $\gamma_0$ and $\beta^{F_1}_{n_1-n_0}(\alpha_0,\dots,\alpha_{n_1-n_0-1})$ is $k$.

 The \textbf{configuration} of $\gamma_0,\gamma_1$ is the pair  $(i',k') \in 2 \times \omega$ given by
 \[
 i' =
 \begin{cases}
 0 &\mbox{ if } A \in \mathscr{U}^{n_1-n_0}\\
 1 &\mbox{ otherwise. }
 \end{cases}
 \]
 and $k'$ is the unique $k < n-2$ for which $B_{F_0,F_1,k} \in \mathscr{U}^{n_1-n_0}$.
\end{enumerate}
\end{definition}

Using the above definition one can prove the following:

\begin{theorem}
\label{thmnonmagidor} Assume AD.
Let $\kappa$ be a strong partition cardinal below $\Theta$, let $(\kappa_n:n\in\omega)$ be the corresponding Kleinberg sequence and let $\kappa_\omega=\bigcup_{n\in\omega}\kappa_n$.
\begin{enumerate}
\item [$(\aleph)$] $\kappa_n$ is not Magidor for every $n>2$.
\item [$(\beth)$] There is a coloring $c:[\kappa_\omega]^{\omega\text{-bd}}\rightarrow\kappa_\omega$ such that if $H\subseteq\kappa_\omega$ and $|H\cap\kappa_n|=\kappa_n$ for some $n>2$ then $c''[H]^{\omega\text{-bd}}=\kappa_\omega$.
\end{enumerate}
\end{theorem}

\par\noindent\emph{Proof}. \newline
Fix $n>2$, and fix a normal ultrafilter over $\mathscr{U}$ over $\kappa$.
We shall define a coloring $c_n:[\kappa_n-\kappa_{n-1}]^{\omega\text{-bd}} \rightarrow{}^\omega 2$.
Given an element $\bar{\nu}=(\nu_\ell:\ell\in\omega)$ in $[\kappa_n]^{\omega\text{-bd}}$ we will define $c_n(\bar{\nu})(m)$ for every $m\in\omega$ thus obtaining an element of ${}^\omega 2$.

Suppose, therefore, that $\nu\in\kappa_n$.
Recall that $\nu=[F]_{\mathscr{U}}$ for some $F:\kappa^{n-2}\rightarrow\kappa$.
Let $\beta^F_{-n+1} : \kappa^{n-2} \to \kappa_2$ be the induced
$(-(n-1))$-relevant function, given in Definition \ref{def:beta},
and $d_n(\nu)$ be the induced $\mathscr{U}$-stable limit value,
\[d_n(\nu) = \inf_{A \in \mathscr{U}^{n-2}}(\beta^F_{-n+1}"A ).\]
The function $d_n(\nu)$ allows us to record whether the interlacing number of two ordinals $\nu_1,\nu_2 \in \kappa_n$ is the maximal possible value $n-2$.

Now for each $m\in\omega$ and each $\bar{\nu}\in[\kappa_n-\kappa_{n-1}]^{\omega\text{-bd}}$ let:

\begin{displaymath}
c_n(\bar{\nu})(m)= \left\{
\begin{array}{ll} 0 & \textrm{if}\ d_n(\nu_{2m})<d_n(\nu_{2m+1})\\
1 & \textrm{if}\ d_n(\nu_{2m})=d_n(\nu_{2m+1})
\end{array} \right.
\end{displaymath}

Suppose that $H\subseteq\kappa_n-\kappa_{n-1}$ and $|H|=\kappa_n$.
If $\beta\in\kappa_2$ then $|\{\nu\in H:d_n(\nu)\leq\beta\}|<\kappa_n$.
We conclude that for every $\beta\in\kappa_2$ there are cobounded sets of $\nu \in H$ with $d_n(\nu) > \beta$.
Also, since the map $\nu \mapsto d_n(\nu)$ on any set of size $\kappa_n$ cannot be injective, we can find $\nu^\beta_0,\nu^\beta_1\in H$ for which $d_n(\nu^\beta_0)=d_n(\nu^\beta_1)>\beta$.
It follows that $c_n''[H]^{\omega\text{-bd}}={}^\omega 2$, as for every $(r_m:m\in\omega)\in{}^\omega 2$ we can pick inductively $\nu_m$ in pairs, so that for every $m\in\omega$ either $d_n(\nu_{2m})<d_n(\nu_{2m+1})$ or $d_n(\nu_{2m})=d_n(\nu_{2m+1})$ according to the desired target value of $r_m$.

Let $\lambda<\Theta$ be a cardinal and let $h:2^\omega\rightarrow\lambda$ be surjective.
The coloring $h\circ c_n: [\kappa_n-\kappa_{n-1}]^{\omega\text{-bd}}\rightarrow\lambda$ will satisfy $h\circ c_n''[H]^{\omega\text{-bd}}=\lambda$ whenever $|H|=\kappa_n$.
In particular, if $\lambda=\kappa_n$ then this coloring exemplifies the non-Magidority of $\kappa_n$ for every $n>2$.

For the second part of the claim we define a coloring $c:[\kappa_\omega]^{\omega\text{-bd}}\rightarrow{}^\omega 2$.
Given $\bar{\nu}\in[\kappa_\omega]^{\omega\text{-bd}}$ fix $n\in\omega$ such that $\{\nu_n:n\in\omega\}\subseteq^*\kappa_n-\kappa_{n-1}$.
Let $\bar{\eta}=(\eta_n:n\in\omega)$ be the part of $\bar{\nu}$ contained in $\kappa_n-\kappa_{n-1}$, re-enumerated as designated.
Define $c(\bar{\nu})=c_n(\bar{\eta})$.

Suppose that $A\subseteq\kappa_\omega$ and $|A|=\kappa_\omega$.
Fix $n\in\omega$ such that $|H|=\kappa_n$ where $H=A\cap(\kappa_n-\kappa_{n-1})$.
By the above definition $c''[H]^{\omega\text{-bd}}=c_n''[H]^{\omega\text{-bd}}=2^\omega$.
Let $h:2^\omega\rightarrow\kappa_\omega$ be surjective.
Now $h\circ c:[\kappa_\omega]^{\omega\text{-bd}}\rightarrow\kappa_\omega$ exemplifies the required statement.

\hfill \qedref{thmnonmagidor}

Our next goal is to prove that $\kappa_\omega$ is Magidor (in the case of $\aleph_\omega$ we shall prove just $\omega$-Magidority).
We shall need a lemma which generalizes a statement from the proof of \cite[Theorem 4.9]{MR0479903}.

\begin{lemma}
\label{lemkl49} Let $\kappa>\aleph_0$ be a strong partition cardinal and let $\kappa_2={}^\kappa\kappa/\mathscr{U}$ where $\mathscr{U}$ is a normal ultrafilter over $\kappa$.
Suppose that $(c_\alpha:\alpha\in\kappa)$ is a sequence of colorings where $c_\alpha:[\kappa_2]^{\kappa^{n_\alpha}}\rightarrow\kappa$ and $n_\alpha\in\omega$ for each $\alpha\in\kappa$.
Then one can find $H\in[\kappa_2]^{\kappa_2}$ so that $H$ is $c_\alpha$-monochromatic for every $\alpha\in\kappa$ simultaneously.
\end{lemma}

\par\noindent\emph{Proof}. \newline
We may assume, without loss of generality, that the domain of each $c_\alpha$ is $[\kappa_2]^{\kappa^\omega}$ since $n_\alpha\in\omega$ and by monotonicity.
If $a\in[\kappa_2]^{\kappa^\omega}$ then $(c_\alpha(a):\alpha\in\kappa)$ is an element of ${}^\kappa\kappa$ (or a subset of $\kappa$).
Therefore, one can think of the sequence $(c_\alpha:\alpha\in\kappa)$ as one single coloring $c:[\kappa_2]^{\kappa^\omega}\rightarrow{}^\kappa\kappa$ defined simply by $c(a)=(c_\alpha(a):\alpha\in\kappa)$.
The proof of the lemma will be accomplished once we show that there exists a $c$-monochromatic $H\in[\kappa_2]^{\kappa_2}$.

Define another coloring $d:[\kappa_2]^{\kappa^\omega+\kappa^\omega}\rightarrow\kappa+1$ as follows.
Given $y\in[\kappa_2]^{\kappa^\omega+\kappa^\omega}$ let $y=a^\frown b$ where $a,b\in[\kappa_2]^{\kappa^\omega}$ and $a<b$.
Let $d(y)=\kappa$ iff $c(a)=c(b)$, and let $d(y)=\nu$ iff $c(a)\neq c(b)$ and $\nu\in\kappa$ is the first ordinal so that $c_\nu(a)\neq c_\nu(b)$.

Since $\kappa_2\rightarrow(\kappa_2)^\gamma_\kappa$ for each $\gamma\in\kappa_2$, one can find a $d$-monochromatic $H\in[\kappa_2]^{\kappa_2}$.
Since $\kappa_2$ is regular, the constant value of $d$ cannot be some $\nu\in\kappa$ (otherwise $H$ will produce $\kappa_2$ distinct values $c_\nu(a)$ for $\kappa_2$-many $a\in [H]^{\kappa^\omega}$).
Hence $d''[H]^{\kappa^\omega}=\{\kappa\}$.
It follows that $c(a)=c(b)$ whenever $a,b\in[H]^{\kappa^\omega}$ and hence $H$ is $c_\alpha$-monochromatic for every $\alpha\in\kappa$ as required.

\hfill \qedref{lemkl49}

In the following proposition we deal with colorings from $[\kappa_\omega]^{\aleph_0\text{-bd}}$ into $\kappa_1$ and plan to establish the existence of sets of maximal cardinality which omit colors by reducing the given coloring to colorings of $\kappa_2$, and applying the conclusion of the previous lemma.
Once this is shown, it is straightforward to establish the Magidority of $\kappa_\omega$ (namely, dealing with colorings into $\kappa_\omega$).

To move forward, we do have to tackle a choice related issue concerning the usage of functions $f \in {}^{\kappa^n} \kappa_2$ to represent ordinals in $\kappa_{n+2}$.
Consider for example the case of  $\kappa_3$.
According to its definition, this cardinal is the order-type of the $\mathscr{U}^2$-induced well-ordering  on its associated equivalence classes of functions $F\in {}^{\kappa^2} \kappa$ (i.e., on $\kappa^2$-sequences of ordinals of $\kappa$).

An alternative way to present ordinals of $\kappa_3$ is by thinking on $\kappa$-sequences of ordinals of $\kappa_2$.
It is natural to consider arbitrary function $f \in {}^\kappa \kappa_2$ as representatives of ordinals in $\kappa_3$.
This is advantageous as it enables us to invoke directly the partition properties of $\kappa_2$.
However, there is an issue with this approach.
as $\kappa$-sequence of ordinals of $\kappa_2$ are merely $\kappa$-sequence of $\mathscr{U}$ equivalence classes, and replacing such a sequence with a function $F \in {}^{\kappa^2}\kappa$ seems to require making more choice than one is allowed in the context of \textsf{AD}.
Therefore, while it is possible to define a well-ordering of $\mathscr{U}$-induced equivalence classes of functions (sequences) $f \in {}^\kappa \kappa_2$, it is possible that some of these classes do not coincide the ordinals in $\kappa_3$, which are given by classes of functions $F \in {}^{\kappa^2}\kappa$.

This issue is not essential, and is solvable once realizing that both types of classes produces well-orders of the same ordertype, which is $\kappa_3$.
This is because each of the two embeds by an order-preserving and unbounded function in the other. To see the fact that they are mutually cofinal we need two directions.
One direction is given by the map $F \mapsto \beta^F_{-1}$ from ${}^{\kappa^2}\kappa$ to ${}^\kappa \kappa_2$.
For the opposite direction consider a function  $\kappa$-sequence $(\delta_\alpha:\alpha\in\kappa)$ of ordinals of $\kappa_2$.
Recall that $\kappa_2$ is regular and hence $(\delta_\alpha:\alpha\in\kappa)$ is bounded in $\kappa_2$.
Choose an upper bound, say $\tau\in\kappa_2$, and let $f\in{}^\kappa\kappa$ be an element from the equivalence class of $\tau$ as an ordinal in $\kappa_2$.
Consider now the $\kappa^2$-sequence formed by $\kappa$ consecutive copies of $f$.
This gives rise to the desired element of ${}^{\kappa^2}\kappa$, and all we need is to choose one representative of one equivalence class.
This consideration generalizes to all cardinals $\kappa_n$, $n \geq 3$ and is summarized by the following lemma.

\begin{lemma}
\label{lemrepresentations} For every $n\in\omega$, the well-orderings ${}^{\kappa^{n+1}}\kappa_2$ and ${}^{\kappa^{n+2}}\kappa$ are mutually cofinal and hence can be employed interchangeably.
\end{lemma}

\hfill \qedref{lemrepresentations}

We turn to make several definitions and notations which adjust previous definitions to work with functions $f \in {}^{\kappa^{n+1}}\kappa_2$, replacing the role of functions $F \in {}^{\kappa^{n+2}}\kappa$ in the previous section.
Fix $n\in\omega$ and $\tau\in\kappa$.
Consider a sequence $\bar{\alpha}=(\alpha_0,\ldots,\alpha_{n-1})\in[\kappa]^n$.
Suppose that $f\in{}^{\kappa^n}\kappa_2$ (so $f$ represents an ordinal of $\kappa_{n+2}$).
For every $i\in\tau$ let $f_i(\bar{\alpha})= f(\alpha_0\cdot\tau,\alpha_1,\ldots,\alpha_{n-1})$.
Notice that $f_i$ is the same object like $f$, that is a function from $\kappa^n$ into $\kappa_2$.
In particular, one may consider the equivalence class $[f_i]_{\mathscr{U}}$.

We define a \emph{breaking function} ${\rm bk}^n_\tau: {}^{\kappa^n}\kappa_2\rightarrow{}^\tau({}^{\kappa^n}\kappa_2)$ as follows.
Given $f\in{}^{\kappa^n}\kappa_2$ let ${\rm bk}^n_\tau(f) =([f_i]_{\mathscr{U}}:i\in\tau)$.
In order to obtain Magidority we shall need only $\tau\in\omega_1$ and this will take care of every possible order type of an $\aleph_0$-bounded subset of $\kappa_\omega$.
We fix, therefore, an enumeration $(\bar{\tau}_\nu:\nu\in\omega_1)$ of ${}^{<\omega}\omega_1$.
Stipulate $\bar{\tau}_\nu=(\tau_\nu(n):n<n_\nu)$, so $n_\nu\in\omega$.

\begin{proposition}
\label{propmag} Suppose that:
\begin{enumerate}
\item [$(a)$] $\kappa>\aleph_0$ is a strong partition cardinal.
\item [$(b)$] $(\kappa_n:n\in\omega)$ is the corresponding Kleinberg sequence generated from the normal ultrafilter $\mathscr{U}$ over $\kappa$.
\item [$(c)$] $\kappa_\omega=\bigcup_{n\in\omega}\kappa_n$.
\item [$(d)$] $c:[\kappa_\omega]^{\aleph_0\text{-bd}}\rightarrow\kappa$.
\end{enumerate}
Then by letting $c_\nu=c\upharpoonright \prod_{n<n_\nu}[\kappa_{n+2}]^{\tau_\nu(n)}$ for $\nu\in\omega_1$, there exists a set $A\in[\kappa_\omega]^{\kappa_\omega}$ which is $c_\nu$-monochromatic for every $\nu\in\omega_1$ simultaneously.
\end{proposition}

\par\noindent\emph{Proof}. \newline
For every $\nu\in\omega_1$ let $\zeta_\nu=\bigcup_{n<n_\nu}\kappa^n$, so $\zeta_\nu$ is an ordinal in $\kappa_2$.
We define for each $\nu\in\omega_1$ an auxiliary function $d_\nu:{}^{\zeta_\nu}\kappa_2\rightarrow\kappa$ by:
$$
d_\nu(f_0,\ldots,f_{n_\nu-1}) = c_\nu(\star({\rm bk}^n_{\tau_\nu(n)}(f_n): n<n_\nu)).
$$
Let us try to explain what lies behind this definition.
One can think of $\zeta_\nu$ as an ordinal, but it is more convenient to identify the members of $\zeta_\nu$ with finite sequences of length $n_\nu$ where each element $f_n$ in such a sequence is a function from $\kappa^n$ to $\kappa_2$.
Such a function is the input of ${\rm bk}^n_{\tau_\nu(n)}$, and given a sequence $(f_0,\ldots,f_{n_\nu-1})$ we apply ${\rm bk}^n_{\tau_\nu(n)}$ to $f_n$ for each $n<n_\nu$ thus creating long sequences (of length $\tau_\nu(n))$.
Finally, we take their concatenation.
Notice that we make a tacit use of Lemma \ref{lemrepresentations} in this definition.

Bearing in mind that $\kappa^n$-sequences of ordinals from $\kappa_2$ represent ordinals in $\kappa_{n+2}$ we see that $\star({\rm bk}^n_{\tau_\nu(n)}(f_n): n<n_\nu)$ is an element in $\prod_{n<n_\nu}[\kappa_{n+2}]^{\tau_\nu(n)}$ so we can apply $c_\nu$ to this element.
From Lemma \ref{lemkl49} we infer that there exists $H\in[\kappa_2]^{\kappa_2}$ which is $d_\nu$-monochromatic for every $\nu\in\omega_1$.
Let $\beta_\nu$ be the constant value of $d_\nu$.
Our goal is to manufacture from $H$ a set $A\in[\kappa_\omega]^{\kappa_\omega}$ which is $c_\nu$-monochromatic for every $\nu\in\omega_1$.
This will be done piecewise, as we shall create $A_n\subseteq\kappa_{n+2}-\kappa_{n+1}$ and define $A=\bigcup_{n\in\omega}A_n$.

The recurrent idea is that ordinals in $\kappa_{n+2}$ are $\kappa^n$-sequences from $\kappa_2$ (in our case, ordinals from $H\subseteq\kappa_2$).
The information in the functions $c_\nu$ is coded by the functions $d_\nu$ for which we already have a monochromatic set at our disposal.
It remains, therefore, to translate the objects in the domain of each $c_\nu$ into the objects which lie in the domain of the $d_\nu$s.
This we shall do now.

Firstly we pick $\omega$-many sequences $(\bar{s}_n:1\leq n<\omega)$ of elements of $H$.
The first sequence $\bar{s}_1$ consists of the first $\kappa$ elements of $H$, and if $\bar{s}_n$ has been defined then $\bar{s}_{n+1}$ consists of the first $\kappa^{n+1}$ elements of $H$ above $\bigcup\bar{s}_n$.
Rather than enumerating the ordinals of each $\bar{s}_n$ by the usual order we would like to ascribe a different enumeration to the elements of each $\bar{s}_n$.

Since the length of $\bar{s}_n$ is the ordinal $\kappa^n$ we shall write $\bar{s}_n$ as $\{s^n_{\bar{\alpha}}:\bar{\alpha}\in[\kappa]^n\}$ where these elements are ordered by the descending lexicographical ordering.
Explicitly, if $\bar{\alpha},\bar{\beta}\in[\kappa]^n, \bar{\alpha}=(\alpha_0,\ldots,\alpha_{n-1})$ and $\bar{\beta}=(\beta_0,\ldots,\beta_{n-1})$ then $\bar{\alpha}<_S\bar{\beta}$ iff $\ell\leq n-1$ is the maximal index for which $\alpha_\ell\neq\beta_\ell$ and then $\alpha_\ell<\beta_\ell$.
We arrange the enumeration of $\bar{s}_n$ in such a way that $s^n_{\bar{\alpha}}<s^n_{\bar{\beta}}$ iff $\bar{\alpha}<_S\bar{\beta}$.

Secondly, for every $n\in\omega$ we define an operator $f\mapsto_n f^+$ which maps each $f\in{}^{\kappa^n}\kappa$ to $f^+\in{}^{\kappa^n}\kappa_2$.
More specifically, it sends $f$ to $f^+\in{}^{\kappa^n}H$ (recall that $H\subseteq\kappa_2$).
This is done by the following procedure.
Given $\bar{\alpha}\in\kappa^n$ and $f\in{}^{\kappa^n}\kappa$ we define $\bar{\alpha}(f)=(\alpha_0,\ldots,\alpha_{n-1},f(\bar{\alpha}))$ and then we let $f^+(\bar{\alpha})=s^{n+1}_{\bar{\alpha}(f)}$.
Philosophically, if $f\in{}^{\kappa^n}\kappa$ then $[f]_{\mathscr{U}}$ is an ordinal in $\kappa_{n+1}$ and therefore $[f^+]_{\mathscr{U}}$ represents an ordinal in $\kappa_{n+2}$, so this is a sort of a \emph{push-up operator}.

Notice that if $f<_{\mathscr{U}}g$ then $f^+<_{\mathscr{U}}g^+$.
Fix $n\in\omega$ and $\ell<n$.
Let $\bar{\alpha}=(\alpha_0,\ldots,\alpha_{\ell-1})\in[\kappa]^\ell$ and fix $f\in{}^{\kappa^n}\kappa$.
By the definition of $f^+$ we see that:
$$
\sup\{f^+(\bar{\alpha}^\frown\bar{\beta}): \bar{\beta}\in[\kappa-(\alpha_{\ell-1}+1)]^{n-\ell}\} = \sup\{s^{n+1}_{\bar{\alpha}^\frown\bar{\beta}}: \bar{\beta}\in[\kappa-(\alpha_{\ell-1}+1)]^{n-\ell}\}.
$$
But the righthand side does not depend on $f$, hence this supremum depends only on $\bar{\alpha}$ and not on $f$.
From this fact it follows that if $(\gamma_i:i\in\tau)$ is an increasing sequence of ordinals in $H$ then one can find a $<_{\mathscr{U}}$-increasing sequence of functions $(f_i^+:i\in\tau)$ so that for some $g\in{}^{\kappa^n}\kappa_2$ we have ${\rm bk}^n_\tau(g)=(f_i^+:i\in\tau)$.
Explicitly, if $(f_i^+:i\in\tau)$ is given then we define:
$$
g(\alpha\cdot\tau+i,\alpha_1,\ldots,\alpha_{n-1})= s^n_{\alpha,\alpha_1,\ldots,\alpha_{n-2}, f_i(\alpha,\alpha_1,\ldots,\alpha_{n-1})}.
$$
One can verify that ${\rm bk}^n_\tau(g)=([f_i^+]_{\mathscr{U}}:i\in\tau)$.

For every $n\in\omega$ let $A_n=\{[f^+]_{\mathscr{U}}: f\in{}^{\kappa^n}\kappa, \kappa=\bigcup(\kappa_n:n\in\omega)\}$.
It follows that $A_n\subseteq\kappa_{n+2}-\kappa_{n+1}$ (but notice that $|A_n|=\kappa_{n+1}<\kappa_{n+2}$).
Finally, let $A=\bigcup_{n\in\omega}A_n$.
By the definition of $d_\nu$ and the construction of $A$ we see that $A$ is $c_\nu$-monochromatic for every $\nu\in\omega_1$, so the proof is accomplished.

\hfill \qedref{propmag}

Now we can phrase and prove the following:

\begin{theorem}
\label{thmadmag} Suppose that $\kappa>\aleph_0$ is a strong partition cardinal.
\begin{enumerate}
\item [$(\aleph)$] If $\kappa>\aleph_1$ then $\kappa_\omega$ is Magidor.
\item [$(\beth)$] If $\kappa=\aleph_1$ then $\kappa_\omega=\aleph_\omega$ is $\omega$-Magidor.
\end{enumerate}
\end{theorem}

\par\noindent\emph{Proof}. \newline
Suppose that $\kappa>\aleph_1$ and let $c:[\kappa_\omega]^{\aleph_0\text{-bd}}\rightarrow\kappa_\omega$ be a coloring.
Define $d:[\kappa_\omega]^{\aleph_0\text{-bd}}\rightarrow\kappa$ as follows.
Given $y\in[\kappa_\omega]^{\aleph_0\text{-bd}}$ let $d(y)=c(y)$ if $c(y)\in\kappa$ and let $d(y)=0$ otherwise.
Notice that $d$ is expressible as $\bigcup_{\nu\in\omega_1}c_\nu$, so by Proposition \ref{propmag} there is a set $A\in[\kappa_\omega]^{\kappa_\omega}$ for which $c_\nu$ is monochromatic whenever $\nu\in\omega_1$.
This implies that $|d''[A]^{\aleph_0\text{-bd}}|\leq\aleph_1<\kappa_1$.
Choose $\delta\in\kappa-d''[A]^{\aleph_0\text{-bd}}$ and verify that $\delta\notin c''[A]^{\aleph_0\text{-bd}}$, so $\kappa_\omega$ is Magidor.

This argument will not work if $\kappa=\aleph_1$, but a similar (and simpler) reasoning shows that one can get $\omega$-Magidority in this case.
Indeed, we do not have to consider $\omega_1$-many possible order types, so we can get rid of the $\bar{\tau}_\nu$s and only consider $\omega$-sequences for which we have but countably many partitions of bounded sequences scattered along the $\aleph_n$s.
We conclude, therefore, that $\aleph_\omega$ is $\omega$-Magidor.

\hfill \qedref{thmadmag}

We conclude this section with a couple of open problems.
If a Kleinberg sequence is inaugurated from $\kappa=\aleph_2$ under \textsf{AD}, then the next step is at least $\aleph_{\omega+1}$ since $\aleph_n$ is singular for every $n>2$.
The possibility of $\aleph_\omega$ being Magidor is not necessarily connected with Kleinberg sequences.
However, if one wishes to pursue this direction then it is reasonable to start with a universe in which infinitely many $\aleph_n$s are measurable.

\begin{question}
  \label{qalephomega} Is it consistent with $\mathsf{ZF}$ that $\aleph_\omega$ is a Magidor cardinal?
\end{question}

Recall that $\lambda$ is $\beta$-Magidor (for some $\beta\in\lambda$) iff $\lambda\rightarrow[\lambda]^{<\beta\text{-bd}}_\lambda$.
A cardinal $\lambda$ is strongly Magidor iff $\lambda$ is $\beta$-Magidor for every $\beta\in\lambda$.
If $\lambda$ is I1 then $\lambda$ is strongly Magidor as proved in \cite{MR3750266}.
It is easy to see from the above statements that if one begins with a larger $\kappa$ as a strong partition cardinal then the Magidority degree of $\kappa_\omega$ increases.
However, it is not clear whether strong Magidority obtains in this way.

\begin{question}
\label{qstrongmag} Assume \textsf{AD}.
Does there exist a singular strongly Magidor cardinal below $\Theta$?
\end{question}

Under AD, infinitary combinatorics extends to non-wellorderable sets. For example, Holshouser and Jackson, Chan \cite{RJonsson}, and 
Chan Jackson Trang \cite{NonWoCards} show that the non-wellordered cardinalities  $\mathbb{R}$, $\mathbb{R} \times \kappa$, and  $[\omega_1]^\omega$ are 
J\'onsson. It is not clear whether the notion of $\aleph_0$-bd partition properties and Magidoricity translates well to this setting. 

\newpage

\section{Forcing Magidor filters over singular cardinals}
\label{sec:ForcingMagidorFilters}

The method of the previous section do not seem to natural extend to show that Kleinberg limit cardinals $\kappa_\omega$ carry a Magidor filter. Indeed, the fact that each of the $\kappa_n$, $n \geq 3$ is not Magidor implies that $\kappa_\omega$ cannot carry a normal Magidor filter.

Considering the existence of Magidor filters in models closely related to \textsf{AD},
we show that one can force over a model of $\textsf{AD}+V=L(\mathbb{R})$ to obtain a Magidor filter over a singular cardinal $\lambda$.
If $\lambda<\Theta$ then our forcing will probably destroy \textsf{AD}, but if one assumes the existence of a strong partition cardinal above $\Theta$ then one can force a Magidor filter over a singular cardinal under \textsf{AD}.
To the understanding of the authors, it is not known whether the existence of strong partition cardinal above $\Theta$ in a model of \textsf{AD} is consistent from a known large cardinal assumption.\footnote{The assumption that there exists a strong partition cardinal above $\Theta$ appears in the literature, see e.g. \cite{MR3633802}.}

We shall see below how to force such filters over models of \textsf{AD}.
Remark that the cardinals upon which we force Magidor filters are not of the form $\kappa_\omega$.
Rather, we choose some strong partition cardinal and singularize it.

We begin with some background.
Prikry proved in \cite{MR0262075} that if one forces with Prikry forcing into a measurable cardinal $\kappa$ then one preserves the fact that $\kappa$ is J\'onsson. Moreover, the normal ultrafilter which serves for the pure part of the forcing conditions generates a J\'onsson filter in the generic extension.
Hence despite the fact that $\kappa$ is not measurable after Prikry forcing, it still secures some combinatorial properties which are quite close to measurability.

We wish to do something similar with Magidor filters.
Namely, we would like to begin with a Magidor ultrafilter over a measurable cardinal, to force Prikry (or Magidor) into it and to show that it generates a Magidor filter in the generic extension.
Our first step, therefore, is to show that there are many Magidor filters over regular (in fact, measurable) cardinals under ${\rm \textsf{AD}}+V=L(\mathbb{R})$.

We shall use, again, the order-type version of Magidority, that is $\lambda\rightarrow[\lambda]^{\omega\text{-bd}}_\lambda$.
This means that the colorings are defined on bounded subsets of $\lambda$ \emph{of order type} $\omega$.
Recall that $\kappa$ a \emph{strong partition cardinal} iff $\kappa\rightarrow(\kappa)^\kappa$ and a \emph{weak partition cardinal} iff $\kappa\rightarrow(\kappa)^\gamma$ for every $\gamma\in\kappa$.

\begin{lemma}
\label{lemmfil} Assume ${\rm AD}+V=L(\mathbb{R})$. \newline
There are unboundedly many $\kappa<\Theta$ which carry a $\kappa$-complete (and normal) Magidor ultrafilter.
\end{lemma}

\par\noindent\emph{Proof}. \newline
By \cite{MR611168} there are unboundedly many strong partition cardinals below $\Theta$.
Let us show that if $\kappa<\Theta$ is such a cardinal then $\kappa$ is measurable, and that if $\mathscr{U}$ is a normal ultrafilter over $\kappa$ then $\mathscr{U}$ is Magidor.
Fix $\kappa<\Theta$ such that $\kappa\rightarrow(\kappa)^\kappa$ and hence $\kappa\rightarrow(\kappa)^2$.
From the latter relation it follows that $\kappa$ is regular, and under ${\rm \textsf{AD}}+V=L(\mathbb{R})$ we know that $\kappa$ is measurable and the filter generated by the unbounded $\omega$-closed subsets of $\kappa$ is a normal ultrafilter, see \cite[Theorem 8.27]{MR2768698}.
Denote this ultrafilter by $\mathscr{U}^\kappa_\omega$, or just $\mathscr{U}$ if $\kappa$ is clear from the context.

Since $\kappa$ is also a weak partition cardinal we see, in particular, that $\kappa\rightarrow(\kappa)^{\omega\cdot\gamma}$ for every countable ordinal $\gamma$.
From \cite[Fact 2.31]{MR2768700} we infer that $\kappa\rightarrow_{\rm club}(\kappa)^{\omega\text{-bd}}$.
Since $\mathscr{U}$ is normal it extends the club filter over $\kappa$ and hence $\kappa\rightarrow(\mathscr{U})^{\omega\text{-bd}}$, so $\mathscr{U}$ is a Magidor filter as desired.

\hfill \qedref{lemmfil}

Notice that we get a somewhat stronger property than just Magidority, as we elicit a monochromatic set from the ultrafilter.

\begin{theorem}
\label{thmmagsing} Magidor filters over singular cardinals. \newline
Assume AD and $V=L(\mathbb{R})$.
Then one can force the existence of a singular cardinal which carries a Magidor filter.
\end{theorem}

\par\noindent\emph{Proof}. \newline
Let $\kappa$ be a strong (and weak) partition cardinal below $\Theta$.
Let $\mathscr{U}$ be a $\kappa$-complete ultrafilter over $\kappa$, so $\mathscr{U}$ is Magidor as we proved in Lemma \ref{lemmfil}.
Let $\mathbb{P}$ be Prikry forcing with $\mathscr{U}$ and let $G\subseteq\mathbb{P}$ be generic over the ground model.
Let $\mathscr{F}$ be the filter generated from $\mathscr{U}$ in $V[G]$.
We claim that $\mathscr{F}$ is a Magidor filter.

So assume that $\name{c}:[\kappa]^{\omega\text{-bd}}\rightarrow\gamma$ for some $\gamma\in\kappa$, and let $(s,A)$ be an arbitrary condition which forces this fact.
We define a function $f:[\kappa]^{\omega\text{-bd}}\rightarrow\gamma^\omega$ in the ground model as follows.
Given $t\in[\kappa]^{\omega\text{-bd}}$ let:
$$
f(t)=\{\eta\in\gamma:\exists u,v \exists E, t=u\cup v, (s^\frown u,E)\Vdash\name{c}(v)=\eta\}.
$$
One can think of $f(t)$ as an $\omega$-sequence of ordinals of $\gamma$.
For every $n\in\omega$ let $f_n:[\kappa]^{\omega\text{-bd}}\rightarrow\gamma$ be the $n$th place in the sequence $f(t)$ and choose $B_n\in\mathscr{U}$ so that $f_n\upharpoonright[B_n]^{\omega\text{-bd}}$ is constant.
We employ $AC_\omega$ at this stage, a principle which holds under \textsf{AD}.
Let $B=\bigcap_{n\in\omega}B_n$, so $B\in\mathscr{U}$ and $f\upharpoonright[B]^{\omega\text{-bd}}$ is constant.
Hence $\mathbb{P}$ forces that $\name{c}\upharpoonright[B]^{\omega\text{-bd}}$ is constant as well.

\hfill \qedref{thmmagsing}

The above theorem shows that a singular cardinal with countable cofinality may carry a Magidor filter.
In a seminal work, Magidor showed in \cite{MR0465868} that from relevant assumptions of large cardinals one can change the cofinality of a regular cardinal making it a singular cardinal of uncountable cofinality.
Using Magidor forcing in the version of Henle from \cite{MR722169} one can obtain a Magidor filter over a singular cardinal with uncountable cofinality, essentially by the same proof.

Results by Chan and Jackson \cite{DestructionOfAD} and Ikegami and Trang \cite{PreservationOfAD} show that AD is highly sensitive to nontrivial forcing below $\Theta$, and in particular show that forcing over a model of AD a Prikry/Magidor sequence to a measurable cardinal $\kappa <  \Theta$ will cause AD to fail in the generic extension (see Theorem 4.1 in \cite{PreservationOfAD}). 
Of course, under the assumption that there is a strong partition cardinal above $\Theta$, Prikry forcing into this cardinal will preserve \textsf{AD}, thus we can phrase the following:

\begin{corollary}
\label{corabovebigtheta} Assume AD.
Suppose that $\kappa>\Theta$ and $\kappa$ is a strong partition cardinal.
Then one can force a Magidor filter over $\kappa$ while preserving AD.
\end{corollary}

\hfill \qedref{corabovebigtheta}

We do not know if one can force a Magidor filter with the usual version of bounded subsets of any countable order.
The problem seems to be Lemma \ref{lemmfil}, in which we need more than ${\rm AC}_\omega$ in order to take care of all the countable order-types simultaneously.
We conclude with the following:

\begin{question}
\label{qbelow} Is it consistent that some singular cardinal $\lambda<\Theta$ carries a Magidor filter under \textsf{AD}?
\end{question}

This question is particularly interesting at $\aleph_\omega$.
Tryba proved in \cite{MR877853} that if there is a J\'{o}nsson filter over $\lambda$ then there are unboundedly many J\'{o}nsson cardinals below $\lambda$, from which he inferred that there is no J\'{o}nsson filter over $\aleph_\omega$ in \textsf{ZFC}.
Under \textsf{AD} there is a J\'{o}nsson filter over $\aleph_\omega$ and this is very natural since each $\aleph_n$ is J\'{o}nsson.
We saw, however, that $\aleph_n$ is not Magidor whenever $n>2$.
Hence even if the answer to the above question is positive, it is plausible that there are no Magidor filters over $\aleph_\omega$ under \textsf{AD}.
Of course, one can force strong properties over infinitely many $\aleph_n$s (e.g., measurability as done in \cite{MR2940502}) and in such models of \textsf{ZF} perhaps there is a Magidor filter over $\aleph_\omega$.

We conclude with another interesting question about strong Magidority, in the context of the present section.
The idea of preserving square bracket relations using Prikry type forcing notions invites for the following:

\begin{question}
\label{qstrmagfil} Is it consistent for a singular cardinal $\lambda$ to carry a strong Magidor filter? Is it consistent under \textsf{AD}?
\end{question}

\textbf{Acknowledgments:} We would like to thank Arthur Apter and Nam Trang for valuable discussion and comments on the topic of the paper.

\newpage

\bibliographystyle{alpha}
\bibliography{arlist}

\end{document}